\documentclass{article}
\usepackage{amssymb,amsmath}
\usepackage{anysize}
\usepackage{psfrag}
\newtheorem{thm}{Theorem}[section]
\newtheorem{def.}{Definition}[section]

\newtheorem{prop}{Proposition}[section]

\newtheorem{cl}{Claim}[section]

\numberwithin{table}{section}

\begin{document}
\title{Hyperfinite knots\\ via the CJKLS invariant in the thermodynamic limit}
        \author{Pedro Lopes\\
        Departamento de Matem\'atica\\
        Instituto Superior T\'ecnico\\
        Av. Rovisco Pais\\
        1049-001 Lisboa\\
        Portugal\\
        \texttt{pelopes@math.ist.utl.pt}\\
}
\date{July 8, 2006}
\maketitle

\begin{abstract}
We set forth a definition of hyperfinite knots. Loosely speaking,
these are limits of certain sequences of knots with increasing
crossing number. These limits exist in appropriate closures of
quotient spaces of knots. We give examples of hyperfinite knots.
These examples stem from an application of the Thermodynamic Limit
to the CJKLS invariant of knots.
\end{abstract}

\section{Introduction} \label{sect:intro}

\noindent

In this article we set a forth a definition of hyperfinite knots.
This definition was originally motivated by trying to regard the
CJKLS invariant (\cite{jsCetal, {pLopes}})  as a partition
function of Statistical Mechanics (\cite{Callen}), and extracting
its Thermodynamic Limit (\cite{Baxter}).

\subsection{Hyperfinite knots} \label{subsect:hyper}

\noindent

Assume there is an invariant of knots, $f$, which takes values on
a closed metric space, $\mathbf{M}$, and consider a quotient space
of knots, ${\cal K}\sb{f}$. Two knots belong to the same
equivalence class of this quotient if they have the same value of
the invariant, $f$. In particular, knots which are deformable into
each other lie in the same equivalence class. For any knot $K$,
let $K\sp{\sim}$ denote the equivalence class, in ${\cal
K}\sb{f}$, which contains $K$. The map $f$ induces an embedding of
this quotient space into $\mathbf{M}$. This allows us to regard
${\cal K}\sb{f}$ as a subspace of $\mathbf{M}$. Let
$\overline{{\cal K}\sb{f}}$ denote the topological closure of
${\cal K}\sb{f}$ with respect to the metric in $\mathbf{M}$.

\bigbreak

Consider a sequence of knots, $\bigl( K\sb{n}\bigr)$, which is a
sequence of representatives of distinct equivalence classes of the
indicated quotient space, with increasing crossing number. Assume
further that the sequence, $\bigl( f (K\sb{n}) \bigr) $ converges
in the closed metric space. There is then an element of
$\overline{{\cal K}\sb{f}}$, call it $K\sb{\infty}$, which is the
limit of the sequence $\bigl( K\sp{\sim}\sb{n}\bigr) $, in the
metric of $M$. We state this as
\[
K\sb{\infty}=\lim\sb{n\rightarrow \infty}K\sp{\sim}\sb{n}
\]
We call $K\sb{\infty}$, as obtained above, a {\it hyperfinite
knot}.

Furthermore, note that the limit of the invariant, $\displaystyle
\lim\sb{n\rightarrow \infty}f(K\sb{n})$, is an invariant of the
{\it hyperfinite knot} $\displaystyle
K\sb{\infty}=\lim\sb{n\rightarrow \infty}K\sp{\sim}\sb{n}$,
keeping the notation above.

\bigbreak

In this article we provide instances where the assumptions above
materialize into concrete examples. This is done by way of a new
invariant of knots which is motivated by applying the so-called
Thermodynamic Limit of Exactly Solved Models (\cite{Baxter}) to
the CJKLS invariant (\cite{jsCetal}).

\subsection{The CJKLS invariant and the thermodynamic limit} \label{subsect:cjklsinv}

\noindent

For each specification of labelling quandle $X$, finite abelian
group $A$ of order $N$, and $2$-cocycle $\phi \in Z\sp{2}\sb{Q}(X,
A)$, the CJKLS invariant (\cite{jsCetal, pLopes}) is a knot
invariant. For each knot $K$ this invariant is denoted $Z\sb{X, A,
\phi}(K)$ or simply $Z(K)$, when the choice of $X$, $A$, and
$\phi$ is clear. It can be regarded as a State-sum where the
states are the homomorphisms (colorings) from the Fundamental
Quandle (\cite{dJoyce, {sMatveev}}) of the knot under study to the
labelling quandle, $X$. The CJKLS invariant takes values in the
group algebra over the integers, $\mathbb{Z}[A]$, of the  abelian
group, $A$. Due to the features of this invariant, the integers
involved in the expression of $Z(K)$ are non-negative, for each
knot $K$. This group algebra, $\mathbb{Z}[A]$, embeds in the group
algebra over the reals of the same finite abelian group,
$\mathbb{R}[A]$. The latter is then regarded as a real vector
space, whose dimension equals the order, $N$, of the finite
abelian group, $A$. We identify this vector space with
$\mathbb{R}\sp{N}$ which, equipped with the Euclidean norm, is a
closed metric space. We, thus, regard the CJKLS invariant as
taking values in this closed metric space.

\bigbreak

 We next think of $Z(K)$ as being a vector State-sum or Partition
 Function in the Helmoltz Representation or Canonical Formalism
 (\cite{Callen}). This means that each coordinate of the vector $Z(K)\in
 \mathbb{R}\sp{N}$ is regarded as a State-sum or Partition Function
 in the Helmoltz Representation or Canonical Formalism. We then
 proceed to extract the Free Energy, $F(K)$. This is the vector
 obtained by taking logarithms of each coordinate of $Z(K)$.
 We drop the ``$-kT$'' factor in the expression of the
 Free Energy since it has no meaning in our set up. We also extend
 the logarithm to have the value zero at zero. Finally, we
 divide each coordinate of the Free Energy by the
 crossing number of $K$ to obtain the Free Energy per crossing,
 $f(K)$. This should correspond to the Free Energy per site in
 Exactly Solved Models (\cite{Baxter}). We remark that this $f(K)$
 is an invariant of $K$, for each knot $K$.

\bigbreak

In this article we fix a labelling quandle, $X$, a finite abelian
group, $A$, and a $2$-cocycle, $\phi$, and consider sequences of
alternating knots whose crossing number strictly increases. For
each such sequence, say $\bigl( K\sb{n}\bigr) $, we calculate the
corresponding sequence of the $f$ invariant, $\bigl(
f(K\sb{n})\bigr) $. Mimicking the Thermodynamic Limit Procedure
(\cite{Baxter}), we calculate the limit of $\bigl(
f(K\sb{n})\bigr) $ as the crossing number goes to infinity. We
obtain infinitely many {\it hyperfinite knots} in this way by
showing that for infintely many sequences $\bigl( K\sb{n}\bigr)$,
$\bigl( f(K\sb{n})\bigr)$ converges.

\bigbreak

The {\it hyperfinite knots} thus establish another liaison between Statistical Mechanics and Knot Theory. In this connection we would like to refer also to the works \cite{aw} and \cite{wda} which relate Exactly Solved Models (\cite{Baxter}) and Knot Theory (\cite{lhKauffman}).

\bigbreak

In \cite{Kauffman0}, Kauffman describes an operator which when
iterated on a knot produces an infinite weaving pattern. He calls
the result of this process an {\it infinite knot}. He formalizes
this limit using a category of infinite sequences. We plan to
investigate the relation of Kauffman's {\it infinite knots} to our
{\it hyperfinite knots} in future work.

\bigbreak

\subsection{Ackowledgements}\label{subsect:ackn}

\bigbreak

\noindent

The author acknowledges support by {\em Programa Operacional
``Ci\^{e}ncia, Tecnologia, Inova\c{c}\~{a}o''} (POCTI) of the {\em
Funda\c{c}\~{a}o para a Ci\^{e}ncia e a Tecnologia} (FCT)
cofinanced by the European Community fund FEDER. He also thanks
the staff at IMPA and especially his host, Marcelo Viana, for
hospitality during his stay at this Institution.

\bigbreak

\section{The definition of hyperfinite knots}\label{sect:infknot}

\noindent

In this section we introduce the formal definitions of the objects
we will be dealing with, {\bf knots} and {\bf hyperfinite knots},
and draw the distinction between them.

\begin{def.}[Knot, \cite{lhKauffman}] A {\bf knot} is an embedding of the standard circle,
$S\sp{1}$, into $\mathbb{R}\sp{3}$. The image of the embedding may
have any $($finite$)$ number of components. Thus, we also use the
word {\bf knot} for what is sometimes called {\bf link} in the
literature.
\end{def.}

\begin{def.}[Crossing number of a knot] Given a knot, $K$, its {\bf crossing
number}, $c\sb{K}$, is the least positive integer with the
following property. There is no diagram of $K$ with strictly less
than $c\sb{K}$ crossings.
\end{def.}

\begin{def.} [Hyperfinite knot] Consider a knot invariant which takes values in a
closed metric space. Specifically, consider a map, $f$, from the
set of knots (or one of its subsets), $\cal K$, into a closed
metric space, $\mathbf{M}$, such that knots that are deformable
into each other receive the same value of $f$.

Let ${\cal K}\sb{f}$ be the quotient set of $\cal K$ by the
relation $\sim $:
\[
K \sim K' \quad \overset{ \text{def.} }{\Longleftrightarrow} \quad
f(K)=f(K')
\]
In particular, knots that are deformable into each other lie in
the same equivalence class. Given $K\in {\cal K}$, let
$K\sp{\sim}\in {\cal K}\sb{f}$ be the equivalence class which
contains $K$. Let $f\sp{\sim} $ be the map from ${\cal K}\sb{f}$
to $M$ which sends $K\sp{\sim}$ to $f(K)$. This map is an
embedding of ${\cal K}\sb{f}$ into $M$. We can thus regard ${\cal
K}\sb{f}$ as a subspace of $M$. Let $ \overline{{\cal K}\sb{f}}$
denote the closure of ${\cal K}\sb{f}$ with respect to the metric
in $\mathbf{M}$.

Assume $\bigl( K\sp{\sim}\sb{n}\bigr)$ is a sequence of knots from
${\cal K}\sb{f}$, with sequence of representatives $\bigl( K\sb{n}
\bigr)$ with increasing crossing number, such that for $n\neq n'$,
$f(K\sb{n})\neq f(K\sb{n'})$. Assume further that $\bigl(
f(K\sb{n})\bigr)$ converges in $M$ to, say $fK\sb{\infty}$. Then
there exists an element of $ \overline{{\cal K}\sb{f}}$, call it
$K\sb{\infty}$ which is the limit of $\bigl(
K\sp{\sim}\sb{n}\bigr)$ in $ \overline{{\cal K}\sb{f}}$:
\[
K\sb{\infty}=\lim\sb{n\rightarrow \infty}K\sb{n}\sp{\sim}
\]
We call such limits {\bf hyperfinite knots}.
\end{def.}

\bigbreak

\begin{prop} Keeping the notation above, $fK\sb{\infty}$ is an
invariant of $K\sb{\infty}$.\end{prop}Proof: It follows from the
definition of hyperfinite knot.
 $\hfill \blacksquare$

\bigbreak

Loosely speaking, {\bf hyperfinite knots} are limits of {\bf
knots} of increasing crossing number. Below, we will give an
example of a proper {\bf hyperfinite knot}, $K\sb{\infty}$, in the
sense that $K\sb{\infty}\in \overline{{\cal K}\sb{f}}\setminus
{\cal K}\sb{f}$.

\section{Background on quandles and the CJKLS invariant}\label{sect:backcjkls}

\subsection{Quandles and the Fundamental Quandle of a Knot}\label{subsect:quandle}

\noindent

\begin{def.} A {\bf quandle}, {\bf\cite{dJoyce, sMatveev}}, is a set, $X$, equipped with a
binary operation, $\ast $, such that, for any $a, b, c \in X$
\begin{itemize}
\item $a \ast a = a$ \item there is a unique $x\in X$ such that
$x\ast b = a$ \item $(a\ast b)\ast c = (a\ast c)\ast (b\ast c)$
\end{itemize}
The quandle formed by the set $X$ and the binary operation $\ast $
will be denoted $\bigl( X, \ast   \bigr)$.
\end{def.}

We remark that the second axiom above gives rise to a second
operation on the quandle denoted $\overline{\ast}$. Then,
$a\overline{\ast}b$ is precisely the unique $x$ guaranteed by the
second axiom such that $x\ast b = a$.

The set of Laurent polynomials in a variable $T$, $\mathbb{Z}[T,
T\sp{-1}]$, endowed with the operation
\[
a \ast b := Ta + (1-T)b \quad \text{ for any } a, b \in
\mathbb{Z}[T, T\sp{-1}]
\]
is an example of a quandle with infinite elements. Quotients of
the set of Laurent polynomials by appropriate ideals and endowed
with the analogous binary operation, give rise to the so-called
Alexander quandles. We consider only the Alexander quandle we will
be interested in, in this article. It is called $S\sb{4}$. This is
the set
\[
\mathbb{Z}\sb{2}[T, T\sp{-1}]/(T\sp{2}+T+1)
\]
endowed with the operation
\[
a\ast b = Ta+(1-T)b \quad \text{ for any } a, b \in
\mathbb{Z}\sb{2}[T, T\sp{-1}]/(T\sp{2}+T+1)
\]
in the indicated quotient. As a set, $\mathbb{Z}\sb{2}[T,
T\sp{-1}]/(T\sp{2}+T+1)$ has four elements whose representatives
may be taken to be $0, 1, T, T+1$. We calculate $1\ast (T+1)$:
\[
1\ast (T+1) = T\cdot 1 +(1-T)\cdot (T+1) = T +T+1-T\sp{2}-T =
2T+T\sp{2}+T+1 = 0
\]
in the indicated quotient. Moreover, since, in the indicated
quotient,
\[
1 = T\sp{2}+T = T(T+1)
\]
then
\[
T\sp{-1}=T+1
\]

\bigbreak

 The three defining axioms of quandles are intimately related to the Reidemeister
 moves of Knot Theory (\cite{lhKauffman}).

\begin{figure}[h!]
    \psfrag{x}{\huge $x$}
    \psfrag{y}{\huge $y$}
    \psfrag{z}{\huge $z=x\ast y$}
    \centerline{\scalebox{.50}{\includegraphics{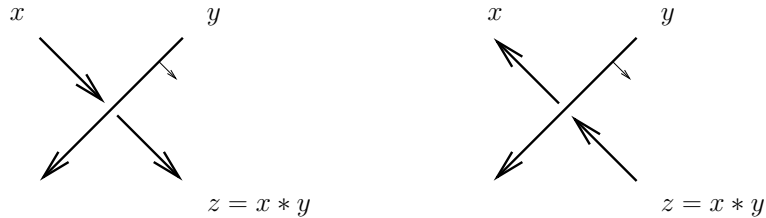}}}
    \caption{Relations read at crossings}\label{Fi:cross0}
\end{figure}

 \begin{def.}The {\bf Fundamental Quandle of a Knot}, $K$,
 $(${\bf\cite{dJoyce, {sMatveev}}}$)$ is presented as follows. Consider any
  diagram of $K$, say $D\sb{K}$, and endow it with an orientation
  and a consistent co-orientation $($i.e., a normal at each point of the diagram$)$.
The arcs of $D\sb{K}$ are regarded as generators and relations of
the following sort are read at each crossing: ``under-arc''
 $\ast $``over-arc'' $ = $ ``under-arc'', where the normal to the over-arc
 points to the under-arc which receives the product, see Figure \ref{Fi:cross0}.
\end{def.}

\begin{thm} The Fundamental Quandle of a knot is a classifying
invariant of it.
\end{thm}Proof: Omitted. See \cite{dJoyce, {sMatveev}}.
$\hfill \blacksquare $

\bigbreak

  The fact that the
 Fundamental Quandle of a Knot $K$ is an invariant of $K$ is a straightforward
 consequence of the defining axioms of the quandle. Unfortunately,
 there is no universal algorithm which allows one to decide, after
 a finite number of steps, whether two presentations of a quandle
 are isomorphic or not. In this way, the Fundamental Quandle of a knot is, per se, of little
 practical use in telling knots apart.

\subsection{Colorings by a labelling quandle}\label{subsect:colorings}

\noindent

Notwithstanding, there is a simple way of obtaining a nice
invariant from the Fundamental Quandle of the Knot. It relies on
the notion of quandle homomorphism.

\begin{def.} Given quandles $\bigl(  X, \ast  \bigr)$ and $\bigl(  X', \ast '
\bigr)$, a homomorphism from $\bigl(  X, \ast  \bigr)$ to $\bigl(
X', \ast ' \bigr)$ is a map, $f$, from $X$ to $X'$ such that, for
any $x, y\in X$,
\[
f(x\ast y)=f(x)\ast ' f(y)
\]
\end{def.}

\begin{def.} Given a knot, $K$, and a finite quandle, $X$, the
{\bf number of colorings of $K$ by $X$} is the number of
homomorphisms from the Fundamental Quandle of the Knot $K$ to $X$.

In this set up, $X$ is referred to as {\bf labelling quandle}. The
elements of the {\bf labelling quandle} are also called {\bf
colors}.
\end{def.}

We remark that there are always the so-called trivial colorings.
These are the colorings where every element of the Fundamental
Quandle of the knot is assigned the same element of the labelling
quandle. There are then at least as many colorings as elements of
the labelling quandle.

\begin{prop} Fix a knot $K$ and a finite quandle $X$. The number
of colorings of $K$ by $X$ is calculated in the following way. We
consider the relations in the presentation of the Fundamental
Quandle of the Knot as a system of equations over the quandle $X$.
The number of solutions of this system of equations is the number
of colorings of $K$ by $X$. Since the Fundamental Quandle of the
Knot is an invariant it then follows that the number of colorings
by $X$ is also an invariant. This number is always at least equal
to the cardinality of the labelling quandle.
\end{prop} Proof: Omitted. See \cite{DL}. $\hfill \blacksquare$

\bigbreak

This invariant was tested in \cite{DL} for its efficiency.

\bigbreak

\subsection{The CJKLS invariant}\label{subsect:cjkls}

\noindent

The CJKLS invariant can be regarded as an elaborated way of
listing the colorings of a knot by a given finite labelling
quandle.

\begin{figure}[h!]
    \psfrag{a}{\huge $a\sb{C}$}
    \psfrag{b}{\huge $b\sb{C}$}
    \psfrag{fi}{\huge $\phi(a\sb{C}, b\sb{C})\sp{+1}$}
    \psfrag{fi'}{\huge $\phi(a\sb{C}, b\sb{C})\sp{-1}$}
    \centerline{\scalebox{.50}{\includegraphics{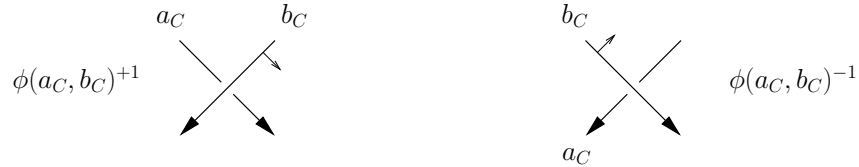}}}
    \caption{The two possible evaluations of $\phi $ at a crossing
($a\sb{C}$ and $b\sb{C}$ are part of an overall coloring
$C$).}\label{Fi:cross}
\end{figure}

\begin{def.} [CJKLS invariant, \cite{jsCetal, {pLopes}}]
\label{def:cjkls1} Choose a finite quandle $X$, a
finite abelian group denoted multiplicatively, $A$, and a
$2$-cocycle $\phi \in Z\sb{Q}\sp{2}(X, A)$ i.e., a map $\phi$ from
$X\times X$ to $A$, such that, for any $a, b, c\in X$
\[
\phi (a, a) = 1 \qquad \qquad \text{ and }  \qquad  \qquad \phi
(a, b)\phi (a\ast b, c)=\phi (a, c)\phi (a\ast c, b\ast c)
\]
where $1$ is the identity in the group $A$.

Given a knot $K$, consider one of its diagrams, $D\sb{K}$, where
the crossings are denoted by $\tau $. Let $\, \cal C$ denote the
set of colorings of the knot $K$ by the labelling quandle $X$.
With respect to the data $X$, $A$, and $\phi $, the CJKLS
invariant of $K$ is
\[
Z\sb{X, A, \phi }(K) = \sum\sb{C\in {\cal C}} \, \prod\sb{\tau \in
D\sb{K}}\phi (a\sb{C}, b\sb{C})\sp{\epsilon\sb{\tau}}
\]
where $\epsilon\sb{\tau}=\pm 1$ and the meaning of $\phi (a\sb{C},
b\sb{C})\sp{\epsilon\sb{\tau}}$ is explained in Figure
\ref{Fi:cross}. When the choice of $X$, $A$, and $\phi$ is clear
we will write $Z(K)$ for $Z\sb{X, A, \phi}(K)$.
\end{def.}

\begin{thm} Keeping the notation above, for each knot $K$, $Z(K)$
is an invariant of $K$.\end{thm} Proof: Omitted. See
\cite{jsCetal, {pLopes}}. $\hfill \blacksquare $

\bigbreak

The fact that, for each knot $K$, $Z(K)$ is an invariant of $K$
stems from the fact that for each coloring $C\in {\cal C}$ the
product
\[
\prod\sb{\tau \in D\sb{K}}\phi (a\sb{C},
b\sb{C})\sp{\epsilon\sb{\tau}}
\]
is invariant under the Reidemeister moves. Unfortunately, since
there is no canonical ordering of the colorings we have to
consider all colorings. Therefore we have to sum the products over
all colorings (\cite{jsCetal}) or simply to list them over all
colorings (\cite{pLopes}).

\bigbreak

We remark that, for any choice of $X$ and $A$, there is always the
so-called trivial $2$-cocycle in $Z\sb{Q}\sp{2}(X, A)$:
\[
\phi\sb{0} \equiv 1
\]
where $1$ is the identity element in the abelian, multiplicative
group $A$. With such a $\phi\sb{0} $ the products in the
expression of $Z$ all become equal to $1$. Then $Z\sb{X, A,
\phi\sb{0}}(K)$ yields the number of colorings of $K$ by the
labelling quandle $X$.

\subsection{An invariant of knots stemming from the CJKLS invariant}\label{subsect:f}

For the remainder of this Subsection, we assume, without loss of
generality, that the following have been fixed. A finite quandle
$X$; a finite abelian group $A$, of order $N$; a $2$-cocycle
$\phi$.

\noindent

\begin{def.}[CJKLS invariant, 2nd version]\label{def:cjkls2}
According to Definition \ref{def:cjkls1}, the CJKLS invariant of
any knot $K$ has the form,
\[
Z(K)=n\sb{1}g\sb{1}+n\sb{2}g\sb{2}+\dots +n\sb{N}g\sb{N}
\]
where the finite abelian group
\[
A=\{ g\sb{1}, g\sb{2}, \dots , g\sb{N}   \}
\]
is ordered in a definite way,  and, for each $i=1, 2, \dots , N$,
\[
n\sb{i} = \# \{ \, C\in {\cal C} \quad |\,  \prod\sb{\tau \in
D\sb{K}}\phi (a\sb{C}, b\sb{C})\sp{\epsilon\sb{\tau}} = g\sb{i} \,
\}\quad (\geq 0)
\]
keeping the notation in Definition \ref{def:cjkls1}.

$Z(K)$ is then an element of $\mathbb{Z}[A]$, the group algebra of
$A$ over the integers. This group algebra embeds in the obvious
way in $\mathbb{R}[A]$ the group algebra of $A$ over the reals.
The latter can be regarded as a real vector space of dimension
$N$, so we identify it with $\mathbb{R}\sp{N}$.

In this way, our second version of the CJKLS invariant is,  for
any knot $K$, to define it in the following way
\[
Z(K):=\bigl( n\sb{1}, n\sb{2}, \dots , n\sb{N} \bigr)
\]
regarded as the (non-negative) coordinates of a vector in
$\mathbb{R}\sp{N}$ with respect to the canonical basis and where
\[
n\sb{i} = \# \{ \, C\in {\cal C} \quad |\,  \prod\sb{\tau \in
D\sb{K}}\phi (a\sb{C}, b\sb{C})\sp{\epsilon\sb{\tau}} = g\sb{i} \,
\}
\]
for each $i=1, 2, \dots , N$, keeping the notation above.
\end{def.}

We will interchangeably use one or the other versions of the CJKLS
invariant.

\bigbreak

\begin{prop}\label{prop:nigt0} We keep the notation above. For any knot $K$, the sum
\[
\sum\sb{i=1}\sp{N}n\sb{i}
\]
equals the number of colorings of the knot by the labelling
quandle. In particular, for at least one $i$
\[
n\sb{i}>0
\]
\end{prop} Proof: Assume
\[
Z(K)=n\sb{1}g\sb{1}+\dots  + n\sb{N}g\sb{N}
\]
If we set each of the $g\sb{i}$'s equal to $1$ this corresponds to
using the trivial $2$-cocycle. But with this choice of
$2$-cocycle, $Z(K)$ yields the number of colorings. The result
follows. $\hfill \blacksquare$

\bigbreak

\begin{def.} We extend the
natural logarithm to zero with value zero and we use the same
symbol ``$\, \ln$'' to denote the extended function.
\end{def.}

\begin{def.}\label{def:Ff}  We keep the notation of Definition
\ref{def:cjkls2}, above.

We define, for any knot $K$,
\[
F(K):=\bigl( \ln (n\sb{1}), \ln (n\sb{2}), \dots , \ln
(n\sb{N})\bigr)
\]
and
\[
f(K):=\biggl( \frac{\ln (n\sb{1})}{c\sb{K}}, \frac{\ln
(n\sb{2})}{c\sb{K}}, \dots , \frac{\ln (n\sb{N})}{c\sb{K}}\biggr)
\]
where $c\sb{K}$ is the crossing number of $K$, and
\[
n\sb{i} = \# \{ \, C\in {\cal C} \quad |\,  \prod\sb{\tau \in
D\sb{K}}\phi (a\sb{C}, b\sb{C})\sp{\epsilon\sb{\tau}} = g\sb{i} \,
\}
\]
for each $i=1, 2, \dots , N$.
\end{def.}

\begin{prop} For each knot $K$, $F(K)$ and $f(K)$ are invariants
of knots which take values in the closed metric space
$\mathbb{R}\sp{N}$, keeping the notation above.
\end{prop} Proof: Omitted. $\hfill \blacksquare$

\bigbreak

We will be particularly interested in the $f$ invariant in the
sequel.

\bigbreak

\subsection{Calculating the CJKLS invariant and the $f$ invariant }\label{subsect:calc}

\noindent

 In \cite{jsCetal0} we find the labelling quandle, $X$, the
abelian group, $A$,  and the $2$-cocycle, $\phi$, we will be
working with in this article.

The labelling quandle, $X$, is the so-called $S\sb{4}$ (page 47)
which is identified with the Alexander quandle
$\mathbb{Z}\sb{2}[T, T\sp{-1}]/(T\sp{2}+T+1)$ (page 48) with
quandle operation
\[
a \ast b := Ta+(1-T)b
\]
in the indicated quotient.

The abelian group, $A$, is $Z\sb{2}\cong (\, t \, | \, t\sp{2}\,
)$ and the $2$-cocycle, $\phi $, is (page 52):
\[
\phi (a, b) := t\sp{\chi\sb{(0, 1)}(a, b)+\chi\sb{(0, T+1)}(a,
b)+\chi\sb{(1, 0)}(a, b)+\chi\sb{(1, T+1)}(a, b)+\chi\sb{(T+1,
0)}(a, b)+\chi\sb{(T+1, 0)}(a, b)}
\]

\begin{def.}\label{def:s4z2}
In the sequel {\bf CJKLS invariant} will mean the CJKLS invariant
with the choice of $X$, $A$, and $\phi$ above.

The invariants $F$ and $f$, introduced in Definition \ref{def:Ff}
will also refer to this choice of $X$, $A$, and $\phi$.

Moreover, when using the second version of the CJKLS invariant we
will order $A\cong (\, t \, |\, t\sp{2}\, )$ so that $g\sb{1}=1$
and $g\sb{2}=t$. Since the cardinality of this group is $2$, then
$N=2$ and the underlying closed metric space is $\mathbb{R}\sp{2}$
endowed with the Euclidean metric.
\end{def.}

We now evaluate the CJKLS invariant and the $f$ invariant of the
trefoil knot and of its mirror image.

\bigbreak

\begin{figure}[h!]
    \psfrag{a}{\huge $a$}
    \psfrag{b}{\huge $b$}
    \psfrag{a*b}{\huge $Ta +(1-T)b$}
    \psfrag{a'}{\huge $Tb+(1-T)[Ta+(1-T)b]=\cdots  = a$}
    \psfrag{b'}{\huge $T[Ta+(1-T)b]+(1-T)a=\cdots =b$}
    \psfrag{fi}{\huge $\phi(a, b)$}
    \psfrag{fi'}{\huge $\phi(b, Ta+(1-T)b)$}
    \psfrag{fi''}{\huge $\phi(Ta+(1-T)b, a)$}
    \centerline{\scalebox{.450}{\includegraphics{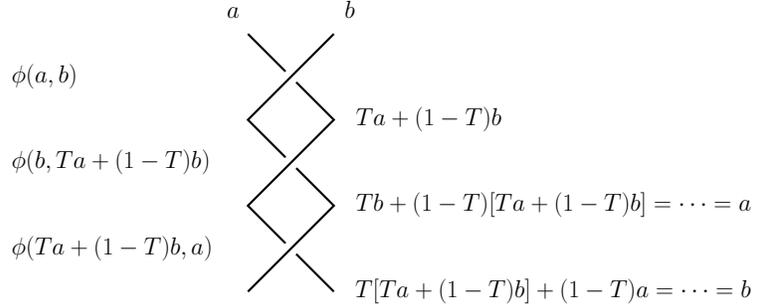}}}
    \caption{The colorings and evaluation of the $2$-cocycle at crossings
for the trefoil}\label{Fi:colortref}
\end{figure}

The trefoil can be regarded as the closure of the braid
$\sigma\sb{1}\sp{3}\in B\sb{2}$ (\cite{Birman}), which is the
braid depicted in Figure \ref{Fi:colortref}. We start by listing
the possible colorings by $S\sb{4}$. In order to do that, we
assign generic colors $a, b\in S\sb{4}$ to the top strands of the
braid in Figure \ref{Fi:colortref} and calculate how they
propagate through each crossing. We enumerate crossings from top
to bottom $1$, $2$, and $3$. The orientation on the strands of the
braid is downwards and the co-orientation is to the left.

The arc emerging from the first crossing is assigned color
\[
a\ast b = Ta+(1-T)b
\]
The arc emerging from the second crossing  is assigned color
\[
b\ast (Ta+(1-T)b) = Tb +(1-T)[Ta+(1-T)b] = (T\sp{2}+T)a +
(T\sp{2}+T+1+2T)b = a
\]
in the indicated quotient. Finally, the arc emerging from the
third crossing is assigned color
\[
(Ta+(1-T)b)\ast a = \dots = b
\]
again, in the indicated quotient. Then, the colors of the strands
at the bottom match the colors of the corresponding strands at the
top, when we close the braid in order to obtain the trefoil. In
this way, any choice of $a$ and $b$ from $S\sb{4}$ gives rise to a
coloring of the trefoil. Since $S\sb{4}$ has four elements, the
number of colorings is $4\sp{2}=16$.

We remark that for a general labelling quandle, what we obtain for
colors at the bottom strands are polynomials in the color inputs
at the top strands. When we equate each of these polynomials to
the corresponding color at the top we obtain a system of equations
in the input colors (\cite{DL}). The number of solutions of this
system of equations is then the number of colorings. In the
present case the system of equations is
\[
\begin{cases}
a=a\\
b=b
\end{cases}
\]

We now calculate, for each $(a, b)\in S\sb{4}\times S\sb{4}$, the
corresponding product of the $\phi $'s over the crossings of the
diagram. The left-hand side of Figure \ref{Fi:colortref} indicates
how to evaluate $\phi$ at each crossing.

In this way,
\[
Z(\text{Trefoil}) = \sum\sb{a, b \in \{ 0, 1, T, 1+T \}} \quad
\phi(a, b)\cdot \phi(b, Ta+(1-T)b)\cdot \phi(Ta+(1-T)b, a)
\]

We recall that the $\phi $ we are using is a function of two
variables from $S\sb{4}$ and takes on values in $Z\sb{2}\cong (\,
t\, |\, t\sp{2}\, )$. Specifically, it takes on $1$ when the two
variables are equal or either one of them equals $T\in S\sb{4}$;
it takes on $t$, otherwise. We then set
\[
\Phi (a, b):=\phi(a, b)\cdot \phi(b, Ta+(1-T)b)\cdot
\phi(Ta+(1-T)b, a)
\]

It is a straightforward exercise to see that
\begin{equation*}
\Phi (a, b)=
\begin{cases}
1, & \text{ if }  a = b \\
t, & \text{ if }   a\neq b
 \end{cases}
\end{equation*}

Thus
\[
\Phi (a, b) = t\sp{\bar{\delta}\sb{a, b}}
\]
with
\begin{equation*}
\bar{\delta}\sb{a, b}=
\begin{cases}
0, & a = b\\
1, & a\neq b
 \end{cases}
\end{equation*}

and so

\[
Z(\text{Trefoil}) = \sum\sb{a, b \in \{ 0, 1, T, 1+T \}}\,
t\sp{\bar{\delta}\sb{a, b}} \quad = \quad 4t\sp{0}+12t \quad =
\quad 4(1 + 3t)
\]

\bigbreak

\begin{figure}[h!]
    \psfrag{a}{\huge $a$}
    \psfrag{b}{\huge $b$}
    \psfrag{bsba}{\huge $(T+1)b+Ta$}
    \psfrag{a'}{\huge $(T+1)a+T[(T+1)b+Ta]=\cdots  = b$}
    \psfrag{b'}{\huge $(T+1)[(T+1)b+Ta]+Tb=\cdots =a$}
    \psfrag{fi}{\huge $\phi(a, b)\sp{-1}$}
    \psfrag{fi'}{\huge $\phi(b, (T+1)b+Ta)\sp{-1}$}
    \psfrag{fi''}{\huge $\phi((T+1)b+Ta, a)\sp{-1}$}
    \centerline{\scalebox{.450}{\includegraphics{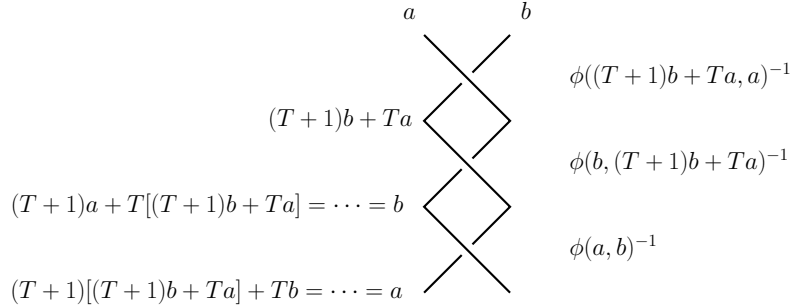}}}
    \caption{The colorings and evaluation of the $\phi $ map at crossings
for the mirror image of the trefoil}\label{Fi:colormtref}
\end{figure}

Now for the mirror image of the trefoil knot. Note that in
$S\sb{4}$
\[
a\overline{\ast}b=T\sp{-1}a+(1-T\sp{-1})b=(T+1)a+Tb
\]
as remarked above. We believe that Figure \ref{Fi:colormtref} is
now self-explanatory.

Then the CJKLS invariant of the mirror image of the trefoil knot
is:
\begin{align*}
\sum\sb{a, b \in \{ 0, 1, T, 1+T \}} &\quad \phi(Ta+(1-T)b,
a)\sp{-1}\cdot \phi(b, Ta+(1-T)b)\sp{-1}\cdot \phi(a, b)\sp{-1}
\quad = \quad \sum\sb{a, b \in \{ 0, 1, T, 1+T \}}\, \Phi (a,
b)\sp{-1} \quad \\
&= \quad \sum\sb{a, b \in \{ 0, 1, T, 1+T \}}\, \Phi (a, b) \quad
= \quad 4(1 + 3t)
\end{align*}
where the equality before the last one follows since $\Phi$ takes
values in $\mathbb{Z}\sb{2}$. Then the CJKLS invariant of both the
trefoil and of its mirror image is $4(1 + 3t)$.

\bigbreak

With our current choice of $X$, $A$ and $\phi$, the CJKLS
invariant is assumed to take values in $\mathbb{R}\sp{2}$, the
first coordinate corresponding to the identity element in
$A=\mathbb{Z}\sb{2}$ and the second coordinate corresponding to
the other element of $A$ (Definition \ref{def:s4z2}). In this way,
\[
f(\text{Trefoil}) = \biggl(  \frac{\ln (4)}{3}, \frac{\ln (12)}{3}
\biggr) = \biggl(  \frac{2\ln (2)}{3}, \frac{2\ln (2) +\ln (3)}{3}
\biggr)
\]
and the same for the mirror image of the trefoil.

\bigbreak

\begin{figure}[h!]
    \psfrag{a}{\huge $a\sb{i-1}$}
    \psfrag{b}{\huge $a\sb{i}$}
    \psfrag{a*b}{\huge $Ta\sb{i-1} +(1-T)a\sb{i}$}
    \psfrag{a*'b}{\huge $T\sp{-1}a\sb{i} +(1-T\sp{-1})a\sb{i-1}$}
    \centerline{\scalebox{.50}{\includegraphics{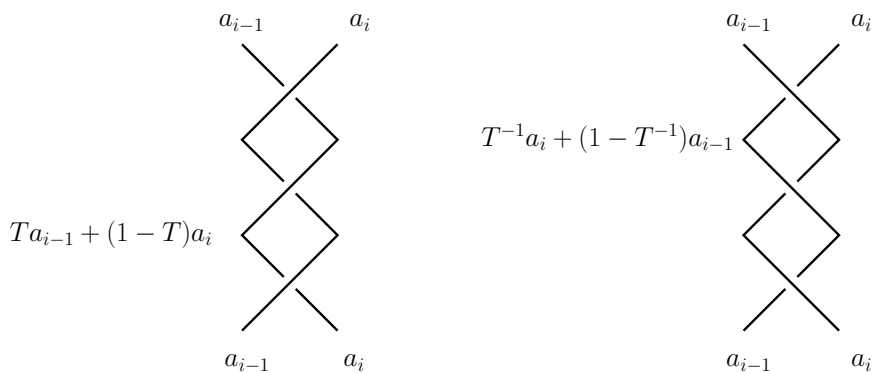}}}
    \caption{The colors of a coloring of $\sigma\sb{i}\sp{\pm 3}$ by $S\sb{4}$}\label{Fi:colortrefi}
\end{figure}

At this point, we record for later use a result whose proof is
implicit in the preceding discussion:

\begin{prop}\label{prop:trefcoloring} Suppose we are calculating
the CJKLS invariant of a given knot using one of its diagrams. In
particular, this diagram has been assigned a coloring by
$S\sb{4}$. Assume further that a certain portion of this diagram
looks like $\sigma\sb{i}\sp{\pm 3}$ and the colors assigned at the
top strands of this $\sigma\sb{i}\sp{\pm 3}$ are $a\sb{i-1},
a\sb{i}$, from left to right, see Figure \ref{Fi:colortrefi}. Then
the colors at the bottom strands of this $\sigma\sb{i}\sp{\pm 3}$
are  $a\sb{i-1}, a\sb{i}$, from left to right. Moreover, the
contribution of the three crossings of this $\sigma\sb{i}\sp{\pm
3}$ in the coloring under study for the summand of the CJKLS
invariant corresponding to this coloring is the factor
\[
\Phi (a\sb{i-1}, a\sb{i})
\]
\end{prop} Proof: Omitted. $\hfill \blacksquare $

\bigbreak

\bigbreak

\section{Sequences of alternating knots} \label{sect:alterknots}

\noindent

In this Section we describe some sequences of alternating knots
(Subsections \ref{subsect:kn}, \ref{subsect:k'n},
\ref{subsect:k0n}, and \ref{subsect:kninfty}), calculate the
corresponding sequences of the $f$ invariant and their limits. The
existence of these limits will imply the existence of hyperfinite
knots, as explained in Section \ref{sect:infknot}. Before, we
define alternating knots and some other particulars of alternating
knots which will concern us in the sequel. We remark that we use
the word {\it knot} to mean both one- or multi-component knots.

\subsection{Alternating knots} \label{subsect:alterknots}

\noindent

\begin{def.}[Alternating knot] An {\bf alternating knot} is a knot such
that one of its diagrams possesses the following property.
Travelling along the diagram, starting at a given point and coming
back to it, and recording at each crossing whether it was passed
over or under, an alternating sequence of  ``over''s and
``under''s is obtained. The {\bf trefoil} $($Figure
\ref{Fi:colortref}$)$ is an example of an {\rm alternating knot}.
\end{def.}

\begin{figure}[h!]
    \psfrag{s}{\LARGE $\rm smoothing$}
    \centerline{\scalebox{.50}{\includegraphics{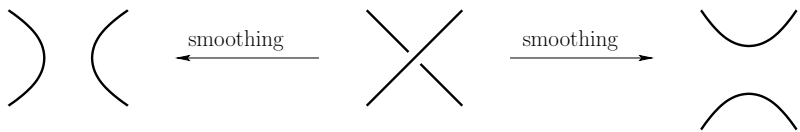}}}
    \caption{The two possible smoothings of a crossing}\label{Fi:crosssmooth}
\end{figure}

\begin{def.}[Smoothing of a crossing] Given a knot diagram, any
crossing of it can be {\bf smoothed} in two different ways, see
Figure  \ref{Fi:crosssmooth}.
\end{def.}

\begin{figure}[h!]
    \psfrag{s}{\huge $smoothing$}
    \centerline{\scalebox{.50}{\includegraphics{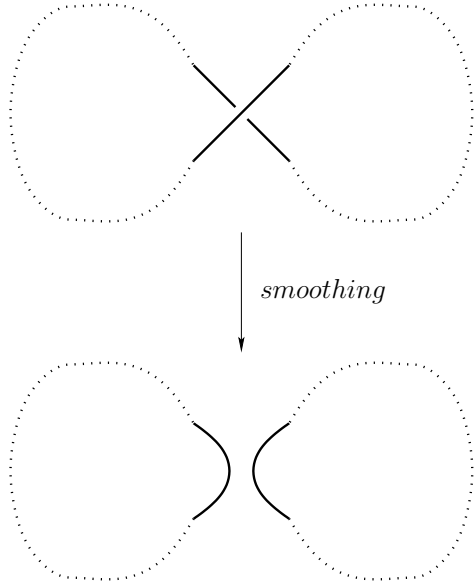}}}
    \caption{A nugatory crossing and the smoothing which reveals its nature}\label{Fi:nugatory}
\end{figure}

\begin{def.}
A {\bf nugatory crossing} of a diagram is a crossing such that one
of its two smoothings disconnects the diagram, see Figure
\ref{Fi:nugatory}. The blank areas surrounded by dotted lines in
this Figure stand for unspecified regions of the diagram.
\end{def.}

\begin{def.}
A {\bf reduced diagram} does not have any {\bf nugatory
crossings}.
\end{def.}

In the sequel, we will need the following property of alternating
knots.

\begin{thm}\label{th:altern} Any two reduced, alternating diagrams of a knot, $K$,
have the same number of crossings. This number is then an
invariant of $K$, the {\bf crossing number of $K$}.
\end{thm}Proof: Omitted. See \cite{Kauffman}, \cite{Murasugi}, or \cite{Thistlethwaite}.
$\hfill \blacksquare $

\bigbreak

We now describe sequences of alternating knots of increasing
crossing number that give rise to hyperfinite knots in the way
described in Section \ref{sect:infknot}.

\bigbreak

\subsection{The $K\sb{n}$ Sequence} \label{subsect:kn}

\noindent

\begin{def.}[The $K\sb{n}$ Sequence] This sequence of alternating knots of increasing crossing number is
given by the closure of the following braids.

\[
b\sb{1} = \sigma\sb{1}\sp{3}, \qquad \sigma\sb{1}\in B\sb{2}
\]
\[
b\sb{2}=\sigma\sb{2}\sp{-3}\cdot \sigma\sb{1}\sp{3}\cdot
\sigma\sb{2}\sp{-3}, \qquad \sigma\sb{1}, \sigma\sb{2}\in B\sb{3}
\]

\[
b\sb{3}=\sigma\sb{3}\sp{3}\cdot \sigma\sb{2}\sp{-3}\cdot
\sigma\sb{1}\sp{3}\cdot \sigma\sb{2}\sp{-3}\cdot
\sigma\sb{3}\sp{3}, \qquad \sigma\sb{1}, \sigma\sb{2},
\sigma\sb{3}\in B\sb{4}
\]

\hspace{200pt}\vdots

\[
b\sb{n}=\sigma\sb{n}\sp{(-1)\sp{n+1}3}\cdot\dots \cdot
\sigma\sb{3}\sp{3}\cdot \sigma\sb{2}\sp{-3}\cdot
\sigma\sb{1}\sp{3}\cdot \sigma\sb{2}\sp{-3}\cdot
\sigma\sb{3}\sp{3}\cdot \dots \cdot
\sigma\sb{n}\sp{(-1)\sp{n+1}3}, \qquad \sigma\sb{1}, \dots ,
\sigma\sb{n}\in B\sb{n+1}
\]
\end{def.}

\begin{figure}[h!]
    \psfrag{a0}{\huge $a\sb{0}$}
    \psfrag{a1}{\huge $a\sb{1}$}
    \psfrag{a2}{\huge $a\sb{2}$}
    \psfrag{a0sa1}{\huge $Ta\sb{0}+(1-T)a\sb{1}$}
    \psfrag{a2sba1}{\huge $T\sp{-1}a\sb{2}+(1-T\sp{-1})a\sb{1}=(T+1)a\sb{2}+Ta\sb{1}$}
    \psfrag{Fi(a1a2)}{\huge $\Phi ( a\sb{1}, a\sb{2})$}
    \psfrag{Fi(a0a1)}{\huge $\Phi ( a\sb{0}, a\sb{1})$}
    \centerline{\scalebox{.450}{\includegraphics{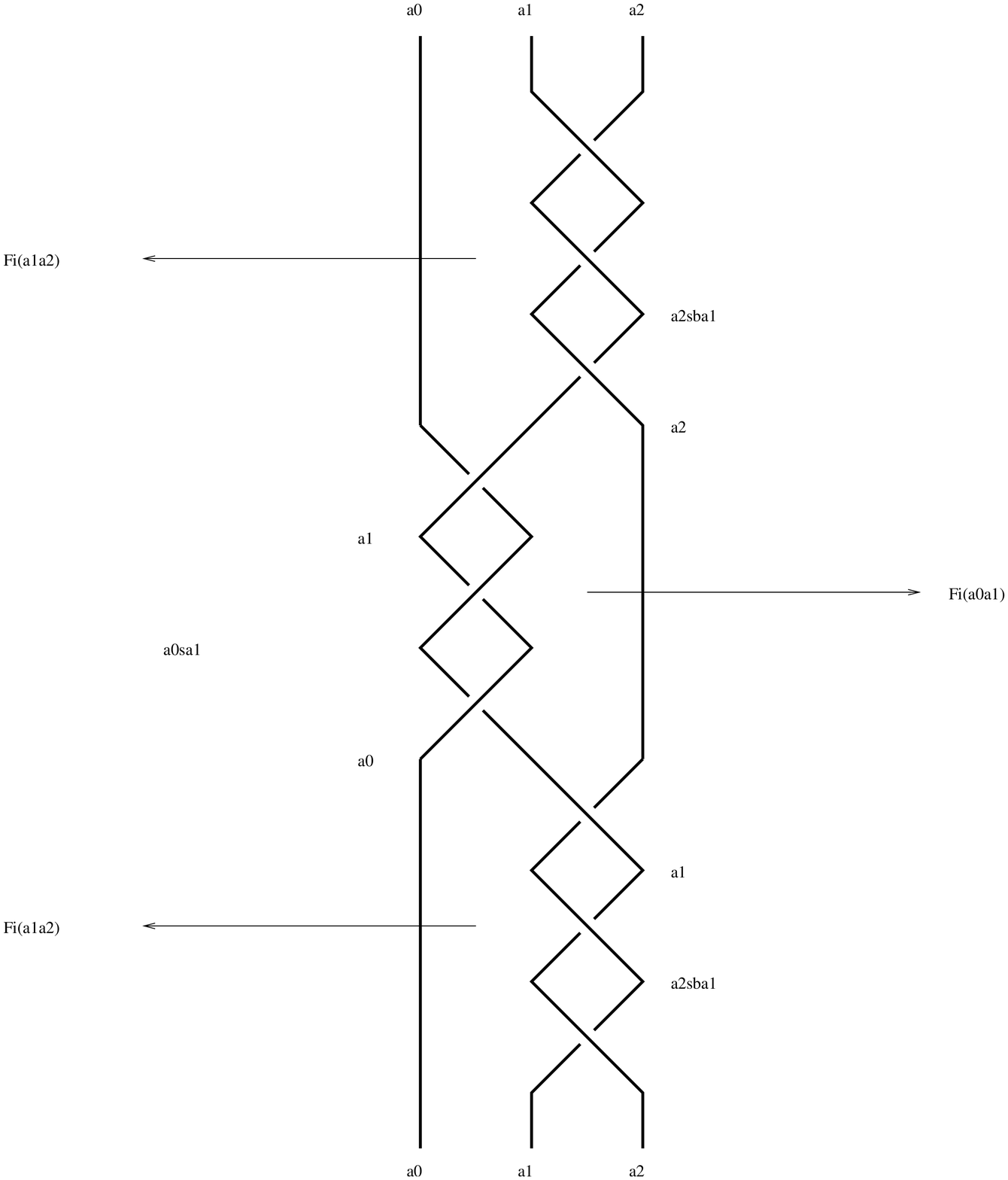}}}
    \caption{$K\sb{2}$, upon closure of the braid, endowed with a coloring by $S\sb{4}$, plus the
$\Phi $ contributions to the CJKLS invariant}\label{Fi:k2}
\end{figure}

\begin{figure}[h!]
    \psfrag{a0}{\huge $a\sb{0}$}
    \psfrag{a1}{\huge $a\sb{1}$}
    \psfrag{a2}{\huge $a\sb{2}$}
    \psfrag{a3}{\huge $a\sb{3}$}
    \psfrag{Fi(a2a3)}{\huge $\Phi ( a\sb{2}, a\sb{3})$}
    \psfrag{Fi(a1a2)}{\huge $\Phi ( a\sb{1}, a\sb{2})$}
    \psfrag{Fi(a0a1)}{\huge $\Phi ( a\sb{0}, a\sb{1})$}
    \centerline{\scalebox{.50}{\includegraphics{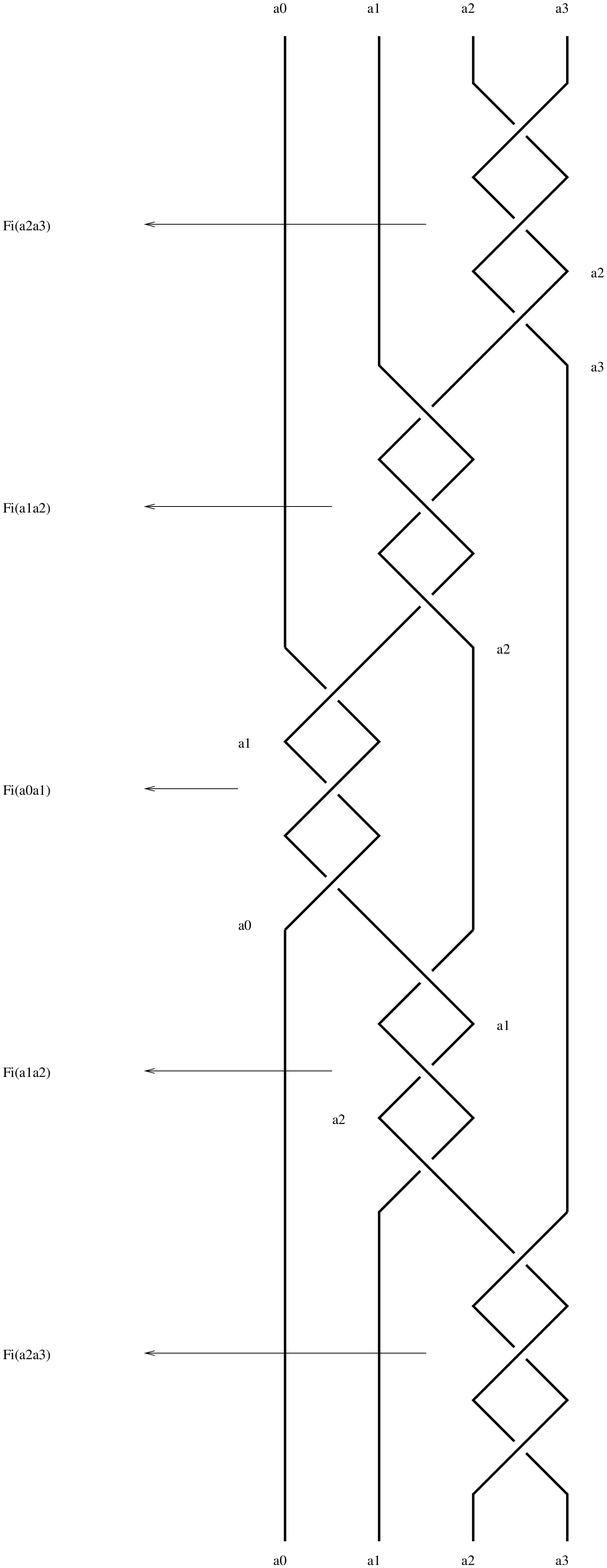}}}
    \caption{$K\sb{3}$, upon closure of the braid,  endowed with a coloring by $S\sb{4}$, plus the
$\Phi $ contributions to the CJKLS invariant}\label{Fi:k3}
\end{figure}

In this sequence, $K\sb{1}$ is the trefoil, regarded as the
closure of the braid $b\sb{1}$ depicted in Figure
\ref{Fi:colortref}. $K\sb{2}$ is the closure of braid $b\sb{2}$
depicted in Figure \ref{Fi:k2}, and $K\sb{3}$ is the closure of
braid $b\sb{3}$ depicted in Figure \ref{Fi:k3}. In this way, the
reader should by now realize how the general $K\sb{n}$ looks like.
We remark that for any $n$, the closure of $b\sb{n}$ is a reduced
alternating diagram. Its number of crossings is then the crossing
number of $K\sb{n}$, according to Theorem \ref{th:altern}.

\bigbreak

As calculated in Section \ref{subsect:calc}, any input $(a, b)$
($a, b\in S\sb{4}$) at the top of the $\sigma\sb{1}\sp{3}$
propagates downwards eventually becoming $(a, b)$ at the bottom of
the braid. In this way there are $4\sp{2}=16$ colorings of the
trefoil by $S\sb{4}$ and the CJKLS invariant is then
\[
Z(\text{Trefoil})=\sum\sb{a, b\in \{ 0, 1, T, T+1
\}}t\sp{\bar{\delta}\sb{a, b}} = 4(1+3t)
\]
where $t$ is the generator of the target group
$\mathbb{Z}\sb{2}\cong (\, t \,  |\, t\sp{2} \, )$. So
\[
Z(K\sb{1})=4+12t \qquad \qquad \qquad f(K\sb{1})=\biggl(
\frac{2\ln (2)}{3}, \frac{2\ln (2) + \ln (3)}{3}\biggr)
\]

The knot $K\sb{2}$ is the closure of the braid $b\sb{2}$ depicted
in Figure \ref{Fi:k2}. Let us first discuss the colorings and the
contribution of each part of the braiding for the CJKLS invariant,
leaning on Proposition \ref{prop:trefcoloring}. The color input is
formed by three colors from $S\sb{4}$, $a\sb{0}, a\sb{1},
a\sb{2}$. The two top right colors $a\sb{1}, a\sb{2}$ propagate
down a $\sigma\sb{2}\sp{-3}$ and so the colors at the bottom of
this $\sigma\sb{2}\sp{-3}$ are $a\sb{1}, a\sb{2}$, according to
Proposition \ref{prop:trefcoloring}. Moreover, the contribution of
this $\sigma\sb{2}\sp{-3}$ to the CJKLS invariant of $K\sb{2}$ is
$\Phi (a\sb{1}, a\sb{2})$, according to the same Proposition. We
now have colors $a\sb{0}, a\sb{1}$ as inputs to
$\sigma\sb{1}\sp{3}$. These colors propagate through
$\sigma\sb{1}\sp{3}$ and so at the bottom of it the colors are
$a\sb{0}, a\sb{1}$. The contribution of this $\sigma\sb{1}\sp{3}$
to the CJKLS invariant of $K\sb{2}$ is $\Phi (a\sb{0}, a\sb{1})$.
Finally, there are colors $a\sb{1}, a\sb{2}$ at the top of
$\sigma\sb{2}\sp{-3}$ and, arguing as before, the contribution of
this $\sigma\sb{2}\sp{-3}$ to the CJKLS invariant of $K\sb{2}$ is
$\Phi (a\sb{1}, a\sb{2})$. The CJKLS invariant of $K\sb{2}$ is
then
\begin{align*}
Z(K\sb{2})&=\sum\sb{a\sb{0}, a\sb{1}, a\sb{2}\in \{ 0, 1, T, T+1
\}}\Phi (a\sb{1}, a\sb{2})\cdot \Phi (a\sb{0}, a\sb{1})\cdot \Phi
(a\sb{1}, a\sb{2}) =\sum\sb{a\sb{0}, a\sb{1}, a\sb{2}\in \{ 0, 1,
T, T+1 \}} \Phi (a\sb{0}, a\sb{1}) = \\
&=\sum\sb{a\sb{0}, a\sb{1}, a\sb{2}\in \{ 0, 1, T, T+1
\}}t\sp{\bar{\delta}\sb{a\sb{0}, a\sb{1}}} =4\, \cdot
\sum\sb{a\sb{0}, a\sb{1}\in \{ 0, 1, T, T+1
\}}t\sp{\bar{\delta}\sb{a\sb{0}, a\sb{1}}} = 4\cdot 4(1+3t) =
4\sp{2}(1+3t)
\end{align*}

\bigbreak

In this way,
\[
Z(K\sb{2}) = 4\sp{2}(1+3t) \qquad \qquad  \qquad f(K\sb{2}) =
\biggl( \frac{4\ln (2)}{9}, \frac{4\ln (2) + \ln (3)}{9}
   \biggr)
\]

The term $K\sb{3}$ is the closure of the braid depicted in Figure
\ref{Fi:k3}. We believe it to be clear now that for any given
integer $n>2$, the contributions from each coloring to the CJKLS
invariant of $K\sb{n}$ come from the $\Phi (a\sb{0}, a\sb{1})$
associated to the $\sigma\sb{1}\sp{3}$, since all other factors
$\Phi (a\sb{i-1}, a\sb{i})$ come in pairs and are thus equal to
$1$ in $\mathbb{Z}\sb{2}$. Then
\[
Z(K\sb{n})= \sum\sb{a\sb{0}, \dots , a\sb{n}\in S\sb{4}}\quad \Phi
(a\sb{0}, a\sb{1}) =\sum\sb{a\sb{0}, \dots , a\sb{n}\in \{ 0, 1,
T, T+1 \}}\quad t\sp{\bar{\delta}\sb{a, b}} \quad = \quad
4\sp{n-1}\cdot 4(1 + 3t) \quad  =  \quad 4\sp{n}(1+3t)
\]

The crossing number  of $K\sb{n}$ is:
\[
3\biggl( \sum\sb{k=1}\sp{n}2  -1\biggr) = 6n-3
\]
We now calculate the $f$ invariant of $K\sb{n}$:
\[
f(K\sb{n})=\biggl(  \frac{\ln (4\sp{n})}{6n-3}, \frac{\ln (3\cdot
4\sp{n})}{6n-3}  \biggr) =\biggl(  \frac{(2n)\ln (2)}{6n-3},
\frac{2n\ln (2) + \ln (3)}{6n-3}  \biggr)
\]

Thus,
\[
Z(K\sb{n})=  4\sp{n}(1+3t)\qquad \qquad \qquad f(K\sb{n})= \biggl(
\frac{(2n)\ln (2)}{6n-3}, \frac{2n\ln (2) + \ln (3)}{6n-3}
\biggr)
\]

\bigbreak

Finally,
\[
\lim\sb{n\rightarrow \infty} f(K\sb{n}) = \lim\sb{n\rightarrow
\infty}\biggl(  \frac{(2n)\ln (2)}{6n-3}, \frac{2n\ln (2) + \ln
(3)}{6n-3}  \biggr) = \biggl( \frac{\ln (2)}{3}, \frac{\ln (2)}{3}
\biggr)
\]

\bigbreak

In this way, the sequence $\bigl( K\sp{\sim}\sb{n}\bigr) $
converges to the hyperfinite knot $K\sb{\infty}$,
\[
K\sb{\infty} = \lim\sb{n\rightarrow \infty}K\sp{\sim}\sb{n}
\]
whose $f$ invariant is
\[
fK\sb{\infty} = \biggl( \frac{\ln (2)}{3}, \frac{\ln (2)}{3}
\biggr)
\]

\bigbreak

\subsection{The $K'\sb{n}$ Sequence} \label{subsect:k'n}

\noindent

\begin{def.}[The $K'\sb{n}$ Sequence] This sequence of alternating knots of increasing crossing number is
given by the closure of the following braids.

\[
b'\sb{1}=b\sb{1} \qquad \qquad b'\sb{2}=b\sb{2} \qquad \qquad
b'\sb{3}=\sigma\sb{3}\sp{3}\sigma\sb{2}\sp{-3}
\sigma{\sb{1}}\sp{3}\sigma\sb{3}\sp{3}\sigma\sb{2}\sp{-3}\sigma\sb{3}\sp{3}
\qquad \qquad
b'\sb{4}=\sigma\sb{4}\sp{-3}\sigma\sb{3}\sp{3}\sigma\sb{2}\sp{-3}
\sigma{\sb{1}}\sp{3}\sigma\sb{3}\sp{3}\sigma\sb{2}\sp{-3}\sigma\sb{3}\sp{3}\sigma\sb{4}\sp{-3}
\]
and in general
\[
b'\sb{2i+1}=\sigma\sb{2i+1}\sp{3}\sigma\sb{2i}\sp{-3}\sigma\sb{2i-1}\sp{3}\cdots
\sigma\sb{3}\sp{3}\sigma\sb{2}\sp{-3}\cdots
\sigma\sb{1}\sp{3}\sigma\sb{3}\sp{3} \cdots
\sigma\sb{2i-1}\sp{3}\sigma\sb{2i+1}\sp{3}\cdots
\sigma\sb{2}\sp{-3}\sigma\sb{3}\sp{3}\cdots
\sigma\sb{2i-1}\sp{3}\sigma\sb{2i}\sp{-3}\sigma\sb{2i}\sp{3}
\]
and
\[
b'\sb{2i+2}=\sigma\sb{2i+2}\sp{-3}b'\sb{2i+1}\sigma\sb{2i+2}\sp{-3}
\]
\end{def.}

\begin{figure}[h!]
    \psfrag{a0}{\huge $a\sb{0}$}
    \psfrag{a1}{\huge $a\sb{1}$}
    \psfrag{a2}{\huge $a\sb{2}$}
    \psfrag{a3}{\huge $a\sb{3}$}
    \psfrag{a4}{\huge $a\sb{4}$}
    \centerline{\scalebox{.50}{\includegraphics{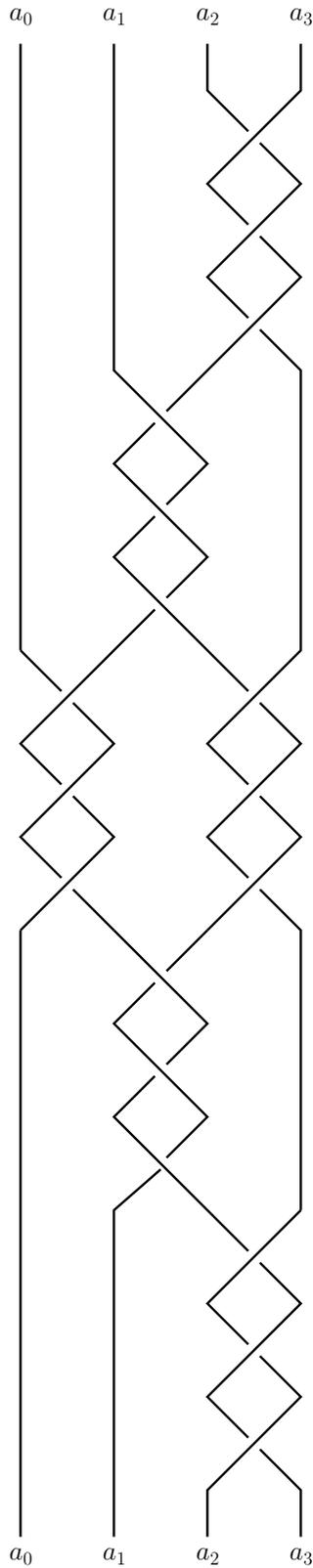}}}
    \caption{$K'\sb{3}$, upon closure of the braid,  endowed with a coloring by $S\sb{4}$}\label{Fi:kprime3}
\end{figure}

\begin{figure}[h!]
    \psfrag{a0}{\huge $a\sb{0}$}
    \psfrag{a1}{\huge $a\sb{1}$}
    \psfrag{a2}{\huge $a\sb{2}$}
    \psfrag{a3}{\huge $a\sb{3}$}
    \psfrag{a4}{\huge $a\sb{4}$}
    \psfrag{a5}{\huge $a\sb{5}$}
    \centerline{\scalebox{.4}{\includegraphics{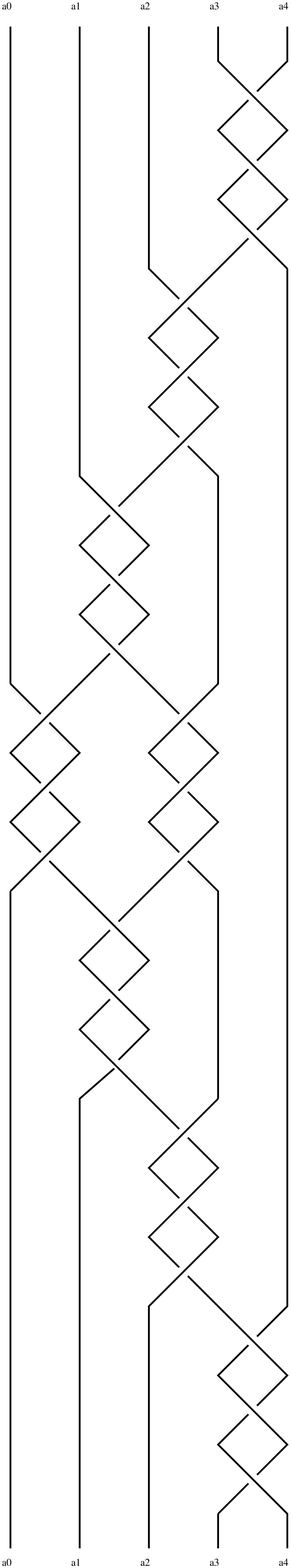}}}
    \caption{$K'\sb{4}$, upon closure of the braid,  endowed with a coloring by $S\sb{4}$}\label{Fi:kprime4}
\end{figure}

\begin{figure}[h!]
    \psfrag{a0}{\huge $a\sb{0}$}
    \psfrag{a1}{\huge $a\sb{1}$}
    \psfrag{a2}{\huge $a\sb{2}$}
    \psfrag{a3}{\huge $a\sb{3}$}
    \psfrag{a4}{\huge $a\sb{4}$}
    \psfrag{a5}{\huge $a\sb{5}$}
    \centerline{\scalebox{.3}{\includegraphics{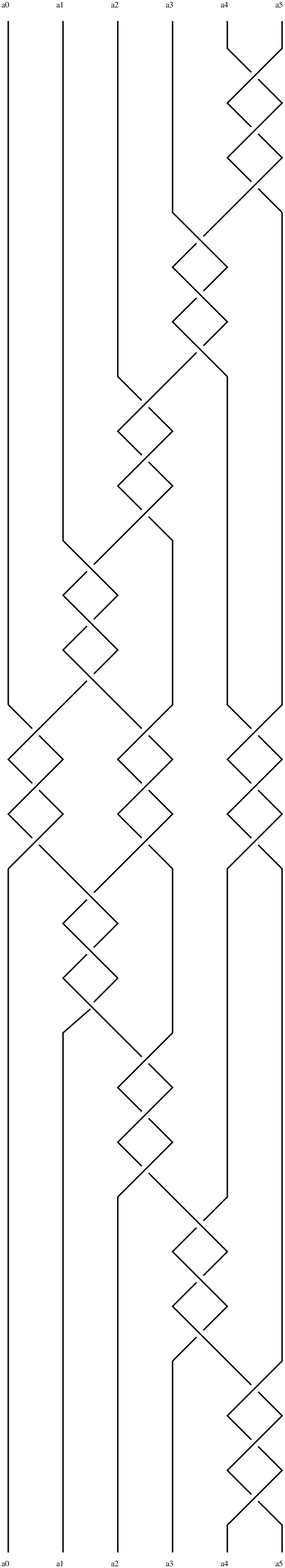}}}
    \caption{$K'\sb{5}$, upon closure of the braid,  endowed with a coloring by $S\sb{4}$}\label{Fi:kprime5}
\end{figure}

We believe the sequence is now clear with the help of Figures
\ref{Fi:kprime3}, \ref{Fi:kprime4}, and \ref{Fi:kprime5} which
depict $b'\sb{3}$, $b'\sb{4}$, and $b'\sb{5}$.

We remark also that for each $n$, the closure of $b'\sb{n}$ is a
reduced alternating diagram. Thus, the number of crossings of
$b'\sb{n}$ is the crossing number of $K'\sb{n}$.

\bigbreak

For odd $n$, the CJKLS invariant is
\begin{align*}
Z(K'\sb{n})&=\sum\sb{a\sb{0}, \dots , a\sb{n}\in \{ 0, 1, T, 1+T
\}}\Phi (a\sb{0}, a\sb{1})\cdot \Phi (a\sb{2}, a\sb{3})\cdot \dots
\cdot \Phi (a\sb{n-1}, a\sb{n}) = \quad \cdots
\end{align*}
where the equality above follows from the fact that the $\Phi $
contributions from the remaining $\sigma\sb{i}\sp{\pm 3}$ come in
pairs and so do not matter.

\begin{align*}
\cdots \quad &=\sum\sb{a\sb{0}, \dots , a\sb{n}\in \{ 0, 1, T, 1+T
\}}t\sp{\bar{\delta}\sb{a\sb{0}, a\sb{1}}}\cdot
t\sp{\bar{\delta}\sb{a\sb{2}, a\sb{3}}}\cdot \dots \cdot
t\sp{\bar{\delta}\sb{a\sb{n-1}, a\sb{n}}} = \quad \cdots
\end{align*}
 We now rewrite this sum in the following way. We write it over
the number, $k$, of pairs $(a\sb{2i}, a\sb{2i+1})$, for $i=0,
\dots , \frac{n-1}{2}$, such that $a\sb{i-1}\neq a\sb{i}$. We
recall that $n$ is odd. We will now count how many possibilities
there are corresponding to $k$ pairs $(a\sb{2i}, a\sb{2i+1})$ with
$a\sb{2i}\neq a\sb{2i+1}$. The set
\[
\{ \,  (a\sb{0}, a\sb{1}), (a\sb{2}, a\sb{3}), (a\sb{4}, a\sb{5}),
\dots , (a\sb{n-1}, a\sb{n})\,  \}
\]
has  $\frac{n+1}{2}$ elements. There are then
\[
\binom{(n+1)/2}{k}
\]
distinct ways of obtaining exactly $k$ pairs $(a\sb{2i},
a\sb{2i+1})$ with distinct coordinates. Since the quandle
$S\sb{4}$ has four elements, for each such pair, $a\sb{2i}$ can
assume one of four elements, whereas $a\sb{2i+1}$ can only assume
one of the remaining three elements. The contribution from $k$
such pairs is $(4\cdot 3)\sp{k}$. Each of the remaining
$\frac{n+1}{2}-k$ pairs has equal coordinates. This can be
realized in four distinct ways for each of them. In this way,
there are
\[
\binom{(n+1)/2}{k}\cdot (4\cdot 3)\sp{k}\cdot 4\sp{(n+1)/2 - k}
\]
distinct ways of realizing exactly $k$ pairs $(a\sb{2i},
a\sb{2i+1})$ with distinct coordinates. Given $k$, the
contribution of the $\Phi$'s to the CJKLS invariant, in each of
these $\displaystyle{\binom{(n+1)/2}{k}\cdot (4\cdot 3)\sp{k}\cdot
4\sp{(n+1)/2 - k}}$ possibilities is $t\sp{k}$ from the pairs with
distinct coordinates and $1\sp{(n+1)/2 - k}$ from the pairs with
equal coordinates. In this way, returning to the evaluation of the
CJKLS invariant:

\begin{align*}
\cdots \quad &=\sum\sb{k=0}\sp{(n+1)/2}\, \binom{(n+1)/2}{k}\cdot
(4\cdot 3)\sp{k}\cdot 4\sp{(n+1)/2 - k}\cdot t\sp{k} =
4\sp{(n+1)/2}\cdot \sum\sb{k=0}\sp{(n+1)/2}\,
\binom{(n+1)/2}{k}\cdot 3\sp{k}\cdot
t\sp{k} =\\
&=4\sp{(n+1)/2}\cdot \sum \sb{k=0, \; \;  k\text{
even}}\sp{(n+1)/2}\, \binom{(n+1)/2}{k}\cdot 3\sp{k} + t\cdot
4\sp{(n+1)/2}\cdot \sum\sb{k=0, \; \;  k\text{ odd}}\sp{(n+1)/2}\,
\binom{(n+1)/2}{k}\cdot
3\sp{k} = \\
&=4\sp{(n+1)/2}\cdot S\sb{(n+1)/2}\sp{e} + t\cdot
4\sp{(n+1)/2}\cdot S\sb{(n+1)/2}\sp{o}
\end{align*}
where
\[
S\sb{m}\sp{e}:=\sum\sb{k=0, \; \;  k\text{ even}}\sp{m}\,
\binom{m}{k}\cdot  3\sp{k}
\]
and
\[
S\sb{m}\sp{o}=\sum\sb{k=0,  \; \;  k\text{ odd}}\sp{m}\,
\binom{m}{k}\cdot 3\sp{k}
\]

\bigbreak

The next result will be useful in the sequel.

\begin{cl}\label{cl:ineq} For integer $m>2$,
\[
3\sp{m}<S\sb{m}\sp{i}<4\sp{m}
\]
for both $i=e$ and $i=o$. \end{cl} Proof: The second inequality is
a consequence of $S\sb{m}\sp{e}, S\sb{m}\sp{o}>0$ and
\[
S\sb{m}\sp{e} + S\sb{m}\sp{o} = (3+1)\sp{m} = 4\sp{m}
\]

As for the first one, let $m>2$ be even.
\begin{align*}
\sum\sb{k=0}\sp{m}\, \binom{m}{k}\cdot 3\sp{k} & = \bigg[
\binom{m}{0}\cdot 3\sp{0}+\binom{m}{2}\cdot 3\sp{2}+\dots +
\binom{m}{m}\cdot 3\sp{m}\bigg] +\bigg[ \binom{m}{1}\cdot 3\sp{1}+
\dots + \binom{m}{m-1}\cdot
3\sp{m-1}\bigg] = \\
& = \bigg[ 1+\frac{m(m-1)}{2}\cdot 3\sp{2}+ \dots + 3\sp{m}\bigg]
+\bigg[ m\cdot 3 + \dots + m\cdot 3\sp{m-1}\bigg]
\end{align*}
Inside each pair of square brackets we find either the summands
corresponding to even $k$ or to odd $k$. Since the last summand
inside the first pair of square brackets is
 $3\sp{m}$ and inside the second pair it is $m\cdot 3\sp{m-1}$ ($m>2$),
the result follows for even $m>2$. Observing that an analogous
calculation holds for odd $m>2$, we conclude the proof. $\hfill
\blacksquare$

\bigbreak

We saw above that, for odd $n$, the CJKLS invariant of $\bigl(
K'\sb{n}\bigr) $ is represented by
\[
\biggl( 4\sp{(n+1)/2}\cdot S\sb{(n+1)/2}\sp{e} \, , \quad
4\sp{(n+1)/2}\cdot S\sb{(n+1)/2}\sp{o}  \biggr)
\]
Since the crossing number of $K'\sb{n}$ is, for odd $n$
\[
3\biggl( \sum\sb{k=1}\sp{n}2  -1\biggr) +3\biggl( \frac{n+1}{2} -
1\biggr) = \frac{12n-6+3n+3-6}{2}=\frac{15n-9}{2}
\]
we then have, for odd $n$,
\[
f(K'\sb{n})=\Biggl( \frac{\ln \bigl( 4\sp{(n+1)/2}\cdot
S\sb{(n+1)/2}\sp{e}\bigr) }{\frac{15n-9}{2}}, \frac{\ln \bigl(
4\sp{(n+1)/2}\cdot S\sb{(n+1)/2}\sp{o}\bigr) }{\frac{15n-9}{2}}
\Biggr)
\]

\bigbreak

For even $n$, we obtain similarly for the CJKLS invariant of
$\bigl( K'\sb{n}\bigr) $,
\[
\biggl( 4\sp{n/2+1}\cdot S\sb{n/2}\sp{e} \, , \quad
4\sp{n/2+1}\cdot S\sb{n/2}\sp{o}  \biggr)
\]
and so, for even $n$
\[
f(K'\sb{n})=\Biggl( \frac{\ln \bigl( 4\sp{n/2+1}\cdot
S\sb{n/2}\sp{e}\bigr) }{\frac{15n-12}{2}}, \frac{\ln \bigl(
4\sp{n/2+1}\cdot S\sb{n/2}\sp{o}\bigr) }{\frac{15n-12}{2}} \Biggr)
\]

\bigbreak

Thanks to Claim \ref{cl:ineq} above, for odd $n$,

\[
0< \frac{\ln
(12)}{15}\underset{n\rightarrow\infty}{\longleftarrow}\frac{\frac{n+1}{2}\ln
(12)}{\frac{15n-9}{2}} \leq \frac{\ln \bigl( 4\sp{(n+1)/2}\cdot
S\sb{(n+1)/2}\sp{i}\bigr) }{\frac{15n-9}{2}}\leq \frac{2(n+1)\ln
(2)}{\frac{15n-9}{2}}\underset{n\rightarrow\infty}{\longrightarrow}\frac{4\ln
(2)}{15}<\frac{\ln (2)}{3}
\]

and analogously, for even $n$

\[
0< \frac{\ln
(12)}{15}\underset{n\rightarrow\infty}{\longleftarrow}\frac{\frac{n+1}{2}\ln
(12) + 2\ln (2)}{\frac{15n-12}{2}} \leq \frac{\ln \bigl(
4\sp{n/2+1}\cdot S\sb{n/2}\sp{i}\bigr) }{\frac{15n-12}{2}}\leq
\frac{2(n+1)\ln
(2)}{\frac{15n-12}{2}}\underset{n\rightarrow\infty}{\longrightarrow}\frac{4\ln
(2)}{15}<\frac{\ln (2)}{3}
\]

We can then conclude that the sequence $\bigl(  f(K'\sb{n}) \bigr)
$  is bounded and so there has to be a convergent subsequence of
it which we denote again by $\bigl(  f(K'\sb{n}) \bigr) $. We can
also conclude that the limit of this convergent sequence, call it
$fK'\sb{\infty}$, is such that
\[
(0, 0) \neq fK'\sb{\infty} \neq fK\sb{\infty}
\]

In this way, the hyperfinite knot
\[
K'\sb{\infty} = \lim\sb{n\rightarrow\infty}K\sp{' \sim}\sb{n}
\]
is different from the hyperfinite knot $K\sb{\infty}$ obtained in
the preceding subsection. Moreover, its $f$ invariant is not $(0,
0)$.

\bigbreak

So far we showed that there are at least two distinct {\it
hyperfinite knots}.

\bigbreak

In the next Subsection we show that the notion of {\it hyperfinite
knot} is non-trivial.

 \bigbreak

\subsection{The $K\sp{0}\sb{n}$ Sequence} \label{subsect:k0n}

\noindent

In this subsection we present a sequence of knots whose $f$
invariant tends to $(0, 0)$.

\begin{def.}[The $K\sp{0}\sb{n}$ Sequence] This sequence of alternating
knots of increasing crossing number is given by the closure of the
braids:
\[
b\sp{0}\sb{1}=\sigma\sb{1}\sp{3} \qquad \sigma\sb{1}\in B\sb{2}
\]
\[
b\sp{0}\sb{2}=\sigma\sb{1}\sp{-3}\sigma\sb{3}\sp{-3}\sigma\sb{2}\sp{3}\sigma\sb{3}\sp{-3}\sigma\sb{1}\sp{-3}
\qquad  \sigma\sb{1}, \sigma\sb{2}, \sigma\sb{3}\in B\sb{4}
\]
\[
b\sp{0}\sb{3}=\sigma\sb{1}\sp{3}\sigma\sb{3}\sp{3}\sigma\sb{5}\sp{3}
\sigma\sb{2}\sp{-3}\sigma\sb{4}\sp{-3}
\sigma\sb{3}\sigma\sb{2}\sp{-3}\sigma\sb{4}\sp{-3}\sigma\sb{1}\sp{3}
\sigma\sb{3}\sp{3}\sigma\sb{5}\sp{3}
\qquad \sigma\sb{1}, \sigma\sb{2}, \sigma\sb{3}, \sigma\sb{4},
\sigma\sb{5}\in B\sb{6}
\]
and in general, for even $n$,
\[
b\sp{0}\sb{n}=\sigma\sb{1}\sp{-3}\sigma\sb{3}\sp{-3}\cdots
\sigma\sb{2n-1}\sp{-3} \sigma\sb{2}\sp{3}\sigma\sb{4}\sp{3}\cdots
\sigma\sb{2n-2}\sp{3}\cdots \sigma\sb{n}\sp{3} \cdots
\sigma\sb{2}\sp{3}\sigma\sb{4}\sp{3}\cdots
\sigma\sb{2n-2}\sp{3}\sigma\sb{1}\sp{-3}\sigma\sb{3}\sp{-3}\cdots
\sigma\sb{2n-1}\sp{-3} \qquad \sigma\sb{1}, \dots ,
\sigma\sb{2n-1}\in B\sb{2n}
\]
and for odd $n$,
\[
b\sp{0}\sb{n}=\sigma\sb{1}\sp{3}\sigma\sb{3}\sp{3}\cdots
\sigma\sb{2n-1}\sp{3} \sigma\sb{2}\sp{-3}\sigma\sb{4}\sp{-3}\cdots
\sigma\sb{2n-2}\sp{-3}\cdots \sigma\sb{n}\sp{3} \cdots
\sigma\sb{2}\sp{-3}\sigma\sb{4}\sp{-3}\cdots
\sigma\sb{2n-2}\sp{-3}\sigma\sb{1}\sp{3}\sigma\sb{3}\sp{3}\cdots
\sigma\sb{2n-1}\sp{3} \qquad \sigma\sb{1}, \dots ,
\sigma\sb{2n-1}\in B\sb{2n}
\]
\end{def.}

\begin{figure}[h!]
    \psfrag{a-2}{\huge $a\sb{-2}$}
    \psfrag{a0}{\huge $a\sb{-1}$}
    \psfrag{a1}{\huge $a\sb{1}$}
    \psfrag{a2}{\huge $a\sb{2}$}
    \psfrag{a3}{\huge $a\sb{3}$}
    \psfrag{a4}{\huge $a\sb{4}$}
    \centerline{\scalebox{.50}{\includegraphics{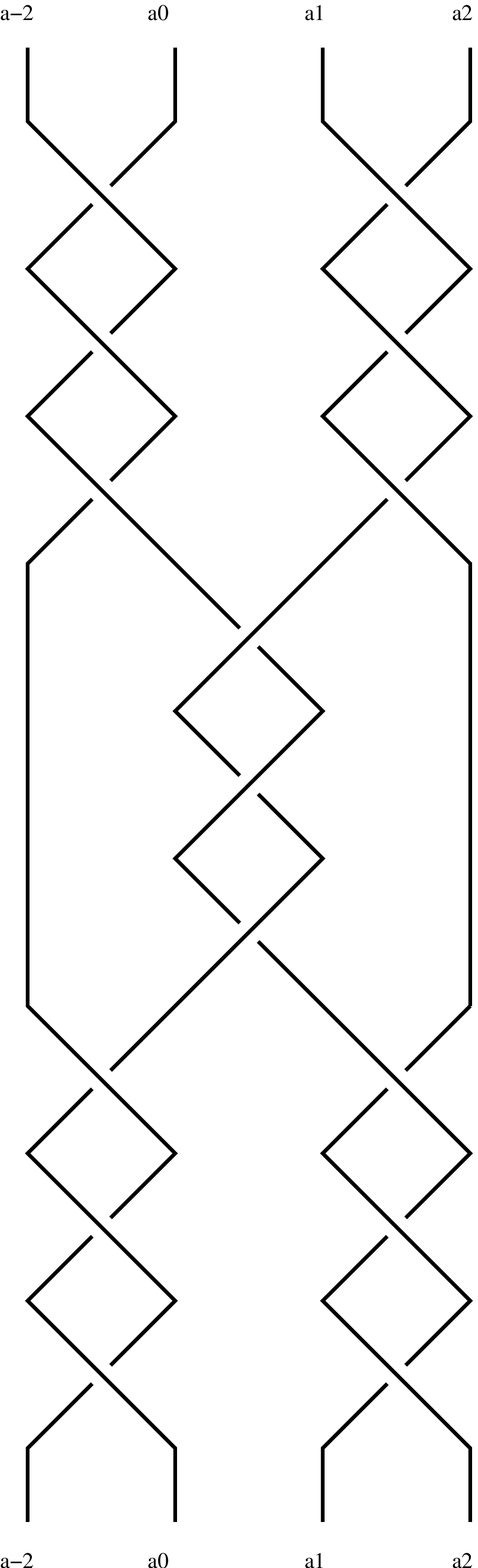}}}
    \caption{$K\sp{0}\sb{2}$, upon closure of the braid,  endowed with a coloring
 by $S\sb{4}$}\label{Fi:k02}
\end{figure}

\begin{figure}[h!]
    \psfrag{a-3}{\huge $a\sb{-3}$}
    \psfrag{a3}{\huge $a\sb{3}$}
    \psfrag{a-2}{\huge $a\sb{-2}$}
    \psfrag{a0}{\huge $a\sb{-1}$}
    \psfrag{a1}{\huge $a\sb{1}$}
    \psfrag{a2}{\huge $a\sb{2}$}
    \psfrag{a3}{\huge $a\sb{3}$}
    \psfrag{a4}{\huge $a\sb{4}$}
    \centerline{\scalebox{.50}{\includegraphics{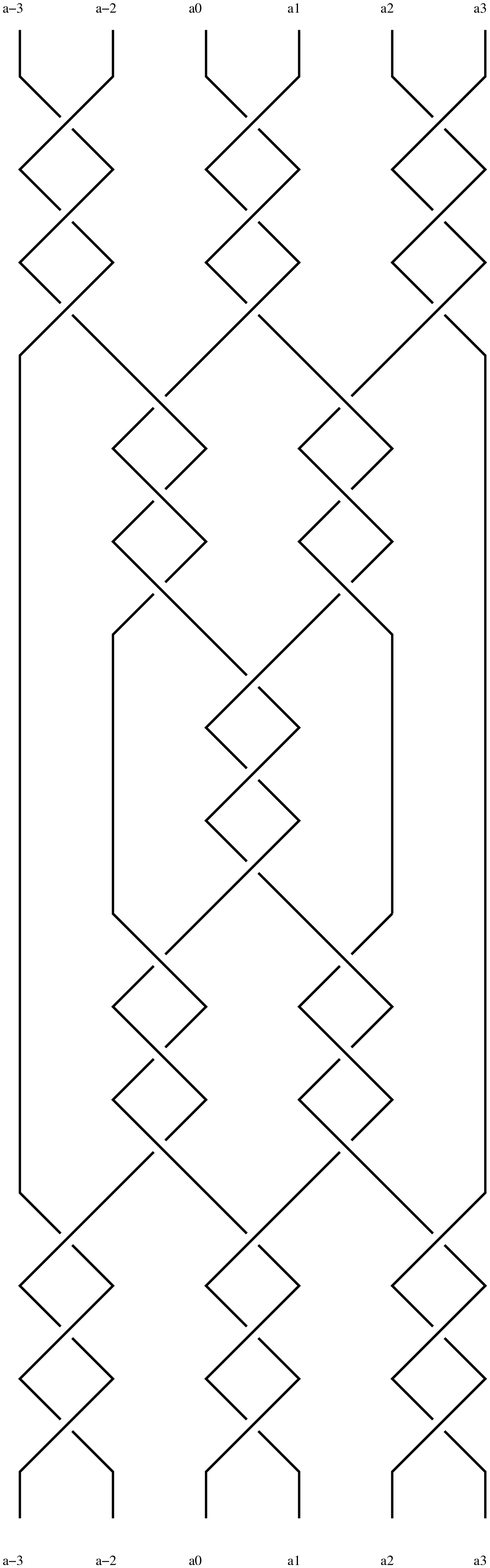}}}
    \caption{$K\sp{0}\sb{3}$, upon closure of the braid,  endowed with a coloring
by $S\sb{4}$}\label{Fi:k03}
\end{figure}

Figures \ref{Fi:k02} and \ref{Fi:k03} should help visualize the
general $K\sp{0}\sb{n}$. The minimum number of crossings of
$K\sp{0}\sb{n}$ is
\[
3\biggl( \sum\sb{k=1}\sp{n}2k \, -1\biggr) = 3\biggl(
2\frac{(n+1)n}{2}\, -1\biggr) = 3n\sp{2}+3n-3
\]

Moreover, using the same sort of analysis as above, the CJKLS
invariant is here
\[
Z\bigl( K\sp{0}\sb{n} \bigr) = 4\sp{2n-2}\cdot 4(1+3t)
\]
thus the $f$ invariant is
\[
f\bigl( K\sp{0}\sb{n} \bigr) = \biggl(  \frac{\ln (4\sp{2n-2}\cdot
4)}{3n\sp{2}+3n-3}, \frac{\ln (4\sp{2n-2}\cdot 4\cdot
3)}{3n\sp{2}+3n-3} \biggr) = \biggl(  \frac{(4n-4)\ln
(4)}{3n\sp{2}+3n-3}, \frac{(4n-4)\ln (2) + \ln (3)}{3n\sp{2}+3n-3}
\biggr)
\]
and so
\[
\lim\sb{n \rightarrow \infty}   f\bigl( K\sp{0}\sb{n} \bigr) =
(0, 0)
\]

There exists then an hyperfinite knot $K\sp{0}\sb{\infty}$ which
is the limit of the sequence $\bigl( K\sp{0 \, \sim}\sb{n}\bigr)$:
\[
K\sp{0}\sb{\infty} = \lim\sb{n\rightarrow \infty} K\sp{0 \,
\sim}\sb{n}
\]

This limit is an example of a {\bf proper hyperfinite knot}.

\begin{prop} $K\sp{0}\sb{\infty}$ belongs
to $\overline{{\cal K}\sb{f}}\setminus {\cal K}\sb{f}$.
\end{prop} Proof: Assume to the contrary an suppose this
hyperfinite knot belongs to ${\cal K}\sb{f}$. Then there should be
a representative of the class with CJKLS invariant equal to
$0+0\cdot t$. But, according to Proposition \ref{prop:nigt0}, at
least one of the $n\sb{i}$'s has to be greater than zero. The
result follows. $\hfill \blacksquare$

 \bigbreak

\subsection{Infinitely many hyperfinite knots} \label{subsect:kninfty}

\noindent

In this subsection we construct two distinct sequences, the
unprimed and the primed sequences of alternating knots, by
replacing the $\sigma\sb{i}\sp{\pm 3}$ by $\sigma\sb{i}\sp{\pm
3(2m+1)}$ in the sequences $\bigl( K\sb{n}\bigr)$ and $\bigl(
K'\sb{n}\bigr)$ of Subsections \ref{subsect:kn} and
\ref{subsect:k'n}, respectively.

\begin{def.} [The $K\sp{m}\sb{n}$ sequence] Let $m$ be a positive integer.
Each term of the sequence $\bigl( K\sp{m}\sb{n}\bigr)$ is given by
the closure of the corresponding term of the sequence of braids:
 \[
b\sp{m}\sb{1}= \sigma\sb{1}\sp{3(2m+1)}, \qquad  \sigma\sb{1}\in
B\sb{1}
\]
\[
b\sp{m}\sb{2}= \sigma\sb{2}\sp{-3(2m+1)}\cdot
\sigma\sb{1}\sp{3(2m+1)}\cdot \sigma\sb{2}\sp{-3(2m+1)}, \qquad
\sigma\sb{1}, \sigma\sb{2} \in B\sb{2}
\]
\[
b\sp{m}\sb{3}=\sigma\sb{3}\sp{3(2m+1)}\cdot
\sigma\sb{2}\sp{-3(2m+1)}\cdot \sigma\sb{1}\sp{3(2m+1)}\cdot
\sigma\sb{2}\sp{-3(2m+1)}\cdot \sigma\sb{3}\sp{3(2m+1)},\qquad
\sigma\sb{1}, \sigma\sb{2}, \sigma\sb{3} \in B\sb{3}
 \]

\hspace{200pt}\vdots

 \begin{align*}
&b\sp{m}\sb{n}=\sigma\sb{n}\sp{(-1)\sp{n+1}3(2m+1)}\cdot\dots
\cdot \sigma\sb{3}\sp{3(2m+1)}\cdot \sigma\sb{2}\sp{-3(2m+1)}\cdot
\sigma\sb{1}\sp{3(2m+1)}\cdot \sigma\sb{2}\sp{-3(2m+1)}\cdot
\sigma\sb{3}\sp{3(2m+1)}\cdot \dots \cdot
\sigma\sb{n}\sp{(-1)\sp{n+1}3(2m+1)}, \\
& \sigma\sb{1}, \dots , \sigma\sb{n}\in B\sb{n+1}
 \end{align*}
\end{def.}

\begin{def.}[The $K\sp{' m}\sb{n}$ Sequence] Let $m$ be a positive
integer. Each term of the sequence $\bigl( K\sp{' m}\sb{n}\bigr)$
is given by the closure of the corresponding term of the sequence
of braids:

\[
b\sp{' m}\sb{1}=b\sp{m}\sb{1} \qquad \qquad
b\sp{' m}\sb{2}=b\sp{m}\sb{2} \qquad \qquad
b\sp{' m}\sb{3}=\sigma\sb{3}\sp{3(2m+1)}\sigma\sb{2}\sp{-3(2m+1)}
\sigma{\sb{1}}\sp{3(2m+1)}\sigma\sb{3}\sp{3(2m+1)}\sigma\sb{2}\sp{-3(2m+1)}\sigma\sb{3}\sp{3(2m+1)}
\]
\[
b\sp{' m}\sb{4}=\sigma\sb{4}\sp{-3(2m+1)}\sigma\sb{3}\sp{3(2m+1)}\sigma\sb{2}\sp{-3(2m+1)}
\sigma{\sb{1}}\sp{3(2m+1)}\sigma\sb{3}\sp{3(2m+1)}\sigma\sb{2}\sp{-3(2m+1)}\sigma\sb{3}\sp{3(2m+1)}
\sigma\sb{4}\sp{-3(2m+1)}
\]
and in general
\begin{multline*}
b\sp{' m}\sb{2i+1}=\sigma\sb{2i+1}\sp{3(2m+1)}\sigma\sb{2i}\sp{-3(2m+1)}\sigma\sb{2i-1}\sp{3(2m+1)}\cdots
\sigma\sb{3}\sp{3(2m+1)}\sigma\sb{2}\sp{-3(2m+1)}\cdots
\sigma\sb{1}\sp{3(2m+1)}\sigma\sb{3}\sp{3(2m+1)} \cdots  \\
 \cdots   \sigma\sb{2i-1}\sp{3(2m+1)}\sigma\sb{2i+1}\sp{3(2m+1)}\cdots
\sigma\sb{2}\sp{-3(2m+1)}\sigma\sb{3}\sp{3(2m+1)}\cdots
\sigma\sb{2i-1}\sp{3(2m+1)}\sigma\sb{2i}\sp{-3(2m+1)}\sigma\sb{2i}\sp{3(2m+1)}
\end{multline*}
and
\[
b\sp{' m}\sb{2i+2}=\sigma\sb{2i+2}\sp{-3(2m+1)}b\sp{' m}\sb{2i+1}\sigma\sb{2i+2}\sp{-3(2m+1)}
\]
\end{def.}

\bigbreak

We now merely state the following results. For any positive
integers $m$ and $n$
\[
Z(K\sp{m}\sb{n}) = Z(K\sb{n}) = 4\sp{n}(1+3t)
\]
and
\begin{equation*}
Z(K\sp{' m}\sb{n}) = Z(K'\sb{n}) =
\begin{cases}4\sp{(n+1)/2}\cdot
S\sb{(n+1)/2}\sp{e} + t\cdot 4\sp{(n+1)/2}\cdot
S\sb{(n+1)/2}\sp{o}, & \qquad \text{  for odd $n$}\\
4\sp{n/2+1}\cdot S\sb{n/2}\sp{e} + t\cdot 4\sp{n/2+1}\cdot
S\sb{n/2}\sp{o}, & \qquad \text{  for even $n$}
\end{cases}
\end{equation*}

The crossing number  of $K\sp{m}\sb{n}$ is
\[
c\sb{K\sp{m}\sb{n}}=3(2m+1)\biggl( \sum\sb{k=1}\sp{n}2  -1\biggr)
= 3(2m+1)(2n-1)
\]
and the crossing number of $K\sp{' m}\sb{n}$ is,
\[
c\sb{K\sp{' m}\sb{n}}=
\begin{cases}
3(2m+1)\biggl( \sum\sb{k=1}\sp{n}2 -1\biggr) +3(2m+1)\biggl(
\displaystyle{\frac{n+1}{2}} - 1\biggr)  , & \, \text{  for odd $n$}\\
3(2m+1)\biggl( \sum\sb{k=1}\sp{n}2 -1\biggr) +3(2m+1)\biggl(
\displaystyle{\frac{n}{2}} - 1\biggr)  , & \, \text{  for even
$n$}
\end{cases}
\]
\[
=
\begin{cases}
\displaystyle{\frac{15n-9}{2}}\, (2m+1)  , & \, \text{  for odd $n$}\\
\displaystyle{\frac{15n-12}{2}}\, (2m+1)  , & \, \text{  for even
$n$}
\end{cases}
\]

then
\[
f(K\sp{m}\sb{n})= \biggl(   \frac{\ln (4\sp{n})}{3(2m+1)(2n-1)},
\frac{\ln (4\sp{n}\cdot 3)}{3(2m+1)(2n-1)} \biggr) = \biggl(
\frac{2n\ln (2)}{3(2m+1)(2n-1)}, \frac{2n\ln (2) + \ln
(3)}{3(2m+1)(2n-1)} \biggr)
\]
with
\[
\lim\sb{n\rightarrow \infty} f(K\sp{m}\sb{n})=  \biggl( \frac{\ln
(2)}{3(2m+1)}, \frac{\ln (2)}{3(2m+1)} \biggr)
\]

So, for each $m$, there exists,
\[
K\sp{m}\sb{\infty}=\lim\sb{n\rightarrow\infty}K\sp{m\, \sim}\sb{n}
\]
which yields an infinite collection of distinct hyperfinite knots.

\bigbreak

Also,

\begin{equation*}
f(K\sp{' m}\sb{n}) =
\begin{cases}\displaystyle{\biggl(
\frac{\ln \bigl( 4\sp{(n+1)/2}\cdot S\sb{(n+1)/2}\sp{e}\bigr)
}{\frac{(15n-9)(2m+1)}{2}} \, , \frac{\ln \bigl(
4\sp{(n+1)/2}\cdot
S\sb{(n+1)/2}\sp{o}\bigr) }{\frac{(15n-9)(2m+1)}{2}}\biggr) }, & \, \text{  for odd $n$}\\
\displaystyle{\biggl( \frac{\ln \bigl( 4\sp{n/2+1}\cdot
S\sb{n/2}\sp{e}\bigr) }{\frac{(15n-12)(2m+1)}{2}}, \frac{\ln
\bigl( 4\sp{n/2+1}\cdot S\sb{n/2}\sp{o}\bigr)
}{\frac{(15n-12)(2m+1)}{2}}\biggr) }\, , & \, \text{ for even $n$}
\end{cases}
\end{equation*}

Arguing as in Subsection \ref{subsect:k'n} using Claim
\ref{cl:ineq} we see that, for each $m$, there is a convergent
subsequence of $\bigl( f(K\sp{' m}\sb{n})\bigr) $ which we denote
again, $\bigl( f(K\sp{' m}\sb{n})\bigr) $. Moreover the limit,
$fK\sp{' m}\sb{\infty}$, of $\bigl( f(K\sp{' m}\sb{n})\bigr) $ is
different from either $(0, 0)$ or $ \biggl( \frac{\ln
(2)}{3(2m+1)}, \frac{\ln (2)}{3(2m+1)} \biggr)$. So, for each $m$,
there is a limit
\[
K\sp{' m}\sb{\infty}=\lim\sb{n\rightarrow\infty}K\sp{' m \,
\sim}\sb{n}
\]
which is neither $K\sp{0}\sb{\infty}$ nor $K\sp{m}\sb{\infty}$.

\bigbreak

\section{Directions for further research}\label{sect:epi}

\noindent

In this article we formalize the notion of {\it hyperfinite knot}.
We only consider here {\it hyperfinite knots} that come from the
CJKLS invariant with the indicated labelling quandle, abelian
group and $2$-cocycle. We would like also to use other data for
the CJKLS invariant and from them construct other {\it hyperfinite
knots}. How do {\it hyperfinite knots} relate for different
choices of the data for the CJKLS invariant? In particular, if a
given sequence of alternating knots with increasing crossing
number converges for a given choice of $X$, $A$, and $\phi$, will
it also converge for a different choice of $X$, $A$, and $\phi$?

Another direction of research would be to look for a different way
of obtaining the $f$ invariant. For instance, could we divide by
the determinant of the knot instead of by its crossing number at
the appropriate step? Or, altogether, find other $f$'s that do not
come from the CJKLS invariant?

Finally, it would be interesting to list the different hyperfinite
knots.

We plan to address these and other questions in future work.

\end{document}